\documentclass[article,dvips,a4paper]{ncc}
\usepackage{amsfonts}
\usepackage{amssymb}
\usepackage{latexcad}
\usepackage{nccbbb}
\usepackage{cite}
\usepackage{bm}

\def\R{{\bbbr}}                                   
\def\sqa{\sqcap\!\!\!\!\sqcup}                    
\def\epr{\phantom{.}\hfill$\sqa$\medskip\par}     
\DeclareMathSymbol{\ell}{\mathord}{letters}{96}   
\def\l{\ell}                                      
\def\G{\Gamma}                                    
\def\g{\gamma}                                    
\def\b{\beta}                                     
\def\a{\alpha}                                    
\def\la{\lambda}                                  
\def\suml {\mathop{\sum}   \limits}               
\def\cdc{,\ldots,}                                
\def\1n{1,\ldots,n}                               
\def\beq{\begin{equation}}                        
\def\eeq{\end{equation}}                          
\def\im{\mathrm{i}}                               
\def\x{ }                                         
\def\_#1{\mathop{\hspace{-0.05ex}^{}_{#1}}}
\def\0L{\text{\Large${\bm 0}$}}
\title{\Large Which Digraphs with Ring Structure are Essentially Cyclic?\vspace{-0.2em}}
\author{\large      {Rafig      Agaev$^1$ and
                     Pavel Chebotarev$^2$}
\vspace{-0.2em}\\
{\normalsize Institute of Control Sciences of the Russian Academy of Sciences}\vspace{-0.4em}\\
{\normalsize 65 Profsoyuznaya Street, Moscow 117997, Russia}}
\date{}

\begin{document}
\maketitle
\footnotetext[1]{E-mail: arpo@ipu.ru} 
\footnotetext[2]{Corresponding author.
E-mail: chv@member.ams.org, upi@ipu.ru; phone: +7-495-334-8869; fax: +7-495-420-2016.}
\setcounter{footnote}{2}

\begin{center}
\begin{minipage}[c]{37.5em}\begin{center} Abstract\end{center}

\small We say that a digraph is essentially cyclic if its Laplacian spectrum is not completely real. The essential cyclicity implies the presence of directed cycles, but not vice versa. The problem of characterizing essential cyclicity in terms of graph topology is difficult and yet unsolved. Its solution is important for some applications of graph theory, including that in decentralized control. In the present paper, this problem is solved with respect to the class of digraphs with ring structure, which models some typical communication networks. It is shown that the digraphs in this class are essentially cyclic, except for certain specified digraphs. The main technical tool we employ is the Chebyshev polynomials of the second kind. A~by-product of this study is a theorem on the zeros of polynomials that differ by one from the products of Chebyshev polynomials of the second kind. We also consider the problem of essential cyclicity for weighted digraphs and enumerate the spanning trees in some digraphs with ring structure.

\end{minipage}
\end{center}

\bigskip
\section{Introduction}
\label{s_intr}

As distinct from the Laplacian eigenvalues of ordinary graphs, those of digraphs need not be real. The problem of characterizing all digraphs that have completely real Laplacian spectra is difficult and yet unsolved. Obviously, two classes of such digraphs are: the acyclic digraphs, whose Laplacian matrices have a triangular form, and symmetric digraphs, whose Laplacian matrices are symmetric and positive semidefinite. It can also be observed that the spectrum of a Laplacian matrix $L$ is completely real whenever for some $\varepsilon>0$, 
$I-\varepsilon L$ is the transition matrix of a reversible Markov chain.

The previous results~\cite{AgaChe05LAA} suggest that non-real eigenvalues with noticeable imaginary parts
characterize\x 
digraphs that have directed cycles not ``extinguished'' 
by counter-directional cycles. On the other hand, there are digraphs that have ``undamped'' cycles and
completely real spectra. In general, it is not an easy problem to distinguish, in terms of graph topology,
digraphs with real Laplacian spectra from those having some non-real Laplacian eigenvalues.
The digraphs of the second type are guaranteed to have cycles, and we call them \emph{essentially\/} cyclic. 
Some preliminary results on essentially cyclic weighted digraphs were presented in~\cite{Che06ILAS,CheAga06ALA}.

The aforementioned problem of distinguishing, in terms of graph topology, digraphs with real and partially non-real
spectra is important for applications. In particular, in the
decentralized control of multi-agent systems \cite{VeermanLafferiere,OlfatiFaxMurray07,RenBeard08,CheAga09ARC}, the
absence of non-real Laplacian eigenvalues of the communication digraph implies that the simplest consensus algorithms
are devoid of oscillations.\x 

A rational 
approach to attacking the difficult problem of characterizing essentially cyclic directed graphs is studying various 
classes of cyclic digraphs. In this paper, we investigate the \emph{digraphs with ring structure}. By such a digraph we mean a digraph whose arc set only contains a collection of arcs forming a Hamiltonian cycle and an arbitrary number of arcs that belong to the inverse Hamiltonian cycle.
Digraphs of this type model a class of typical asymmetric communication networks.
We obtain a necessary and sufficient condition of essential cyclicity for the digraphs with ring structure.
According to this condition, such digraphs are essentially cyclic, except for the digraphs whose inverse Hamiltonian cycle lacks two most distant arcs or one arc or no arc.
This study involves the Chebyshev polynomials of the second kind.
As a by-product we obtain a theorem on the zeros of polynomials that differ by one from the products of Chebyshev polynomials of the second kind.

We also consider weighted digraphs and find that in this case, some conditions of essential cyclicity involve the triangle inequality for the square roots of weight differences.
Finally, we enumerate the converging trees in some digraphs with ring structure.

The paper is organized as follows. After the necessary notation and preliminary results (Section~\ref{s_notation}), in Section~\ref{s_prelim} we present three auxiliary lemmas needed to prove the main results.
In Section~\ref{s_main}, we obtain a necessary and sufficient condition of essential cyclicity for the digraphs with ring
structure, i.\,e., for the digraphs consisting of two opposite Hamiltonian cycles from one of which some arcs can be removed. In Section~\ref{s_weightedv EC}, we investigate the essential cyclicity of the simplest weighted digraphs.
Finally, in Section~\ref{s_concl}, we present explicit formulas for the number of converging trees (in-arborescences) in certain digraphs with ring structure and a direct computation, by means of Chebyshev polynomials, of the 
Laplacian spectrum of the undirected cycle, which is usually proved rather than derived.

\section{Notation and basic results} 
\label{s_notation}

The Laplacian matrix of a digraph $\G$ with vertex set $V(\G)=\{\1n\}$ and arc set $E(\G)$ is the
matrix $L=(\l_{ij})\in\R^{n\times n}$ in which, for $j\ne i$, $\l_{ij}=-1$ whenever $(i,j)\in E(\G)$,
otherwise $\l_{ij}=0$; $\l_{ii}=-\suml_{j\ne i}\l_{ij}$,\, $i,j\in V(\G)$.
If a digraph $\G$ is weighted, i.\,e., each arc $(i,j)\in E(\G)$ has a strictly positive weight $w_{ij}$, then in the
Laplacian matrix $L=L(\G)=(\l_{ij})$, for $j\ne i$, $\l_{ij}=-w_{ij}$ whenever $(i,j)\in E(\G)$,
otherwise $\l_{ij}=0$; $\l_{ii}=-\suml_{j\ne i}\l_{ij}$,\, $i,j\in V(\G)$.
Among the papers concerned with the Laplacian matrices of digraphs, we
mention \cite{FiedlerSedlacek58,CheSha97,CheSha981,CheAga02,AgaChe00,AgaChe01,AgaChe05LAA,Wu,CauVee06Comb} and \cite{Chung05}, where the Laplacian matrix is defined differently. 

We say that a (weighted) digraph $\G$ is \emph{essentially cyclic\/} if its Laplacian spectrum contains non-real
eigenvalues. Evidently, every essentially cyclic digraph has at least one directed cycle. Indeed,
otherwise there exists a numbering of vertices such that the Laplacian matrix has a triangular form,
so the Laplacian spectrum consists of the real diagonal entries of this form.

The \emph{Chebyshev polynomial of the second kind}, $P_n(x),$ \emph{scaled on\/} $]\!-\!2,2[$ is the polynomial of degree $n$
defined by
\beq
\label{200507eq2}
P_n(x)={\sin((n+1)\arccos{x\over 2})\over\sqrt{1-{x^2\over4}}},\quad\text{where}\;\; x\in\,]\!-\!2,2[.
\eeq
Using the auxiliary variable $\varphi\in\,]0,\pi[$ such that $x=2\cos\varphi$ (i.\,e.,
$\varphi=\arccos\frac{x}{2}$) one can rewrite \eqref{200507eq2} in the form
\beq
\label{080207eq1}
P_n(x)={\sin((n+1)\varphi)\over\sin\varphi}.
\eeq

The explicit form of the polynomials $P_n(x)$ is (see, e.g., Theorem~1.12 in \cite{Paszkowski}) 
\beq
\label{080207eq2}
P_n(x)=\suml^{[n/2]}_{i=0} {(-1)}^i\binom{n-i}{i}x^{n-2i}
\eeq
and they satisfy the recurrence 
\beq
\label{230207eq3}
 P_n(x)=x P_{n-1}(x)-P_{n-2}(x)
\eeq
with the initial conditions $P_0(x)\equiv 1$ and $P_1(x)\equiv x$. By \eqref{080207eq1}, the roots of  $P_n(x)$
are:
\beq
\label{210507eq1}
x_k=2\cos{\pi k\over n+1},\quad k=\1n.
\eeq

In particular, if $n=2m-1$ and $m>0$ is integer, then
\beq
\label{e_P2m1}
P_{n}(x)
= \prod_{k=1}^{2m-1}\Big(x  -2\cos  {\pi k\over 2m}\Big)
=x\prod_{k=1}^{ m-1}\Big(x  -2\cos  {\pi k\over 2m}\Big)
                    \Big(x  +2\cos  {\pi k\over 2m}\Big)
=x\prod_{k=1}^{ m-1}\Big(x^2-4\cos^2{\pi k\over 2m}\Big).
\eeq

Consider the tridiagonal matrix $M_n\in\R^{n\times n}:$
\beq
\label{280609eq3}
M_n=\left(%
\begin{array}{rrrrrrr}
 2    &-1    &      &      &      &\0L\\
-1    & 2    &-1    &      &      &   \\
      &\ddots&\ddots&\ddots&      &   \\
      &      &\ddots&\ddots&\ddots&   \\
      &      &      &-1    & 2    &-1 \\
\0L   &      &      &      &-1    & 1 \\
\end{array}%
\right).
\eeq
In particular, $M_1=\big(1\big)$.
Let $Z_n(x)$ be the characteristic polynomial of~$M_n$: $Z_n(x)=\det(xI-M_n)$.
The expansion along the first row of $xI-M_n$ for every $n\ge 2$ provides
\beq
\label{010207eq1}
Z_n(x)=(x-2)Z_{n-1}(x)-Z_{n-2}(x),
\eeq
with the initial conditions $Z_0(x)\equiv1$ and $Z_1(x)\equiv x-1$.
The polynomials $Z_n(x)$ play a 
central role in the subsequent considerations.

\section{Auxiliary lemmas}
\label{s_prelim}

It will be shown in Section~\ref{s_main} that the Laplacian characteristic polynomials of the digraphs with
ring structure differ by one from the products of the polynomials $Z_n(x)$. In this section, we prove
three lemmas. Lemma~\ref{010207prop1} connects $Z_n(x)$ with the Chebyshev polynomials $P_{n}(x)$ of the second kind, Lemma~\ref{t_ZnExpli} provides the explicit form of $Z_n(x)$, and Lemma~\ref{270207le} specifies the roots of the polynomials $Z_n(x)\pm 1$.
\begin{lemma}
\label{010207prop1}
For $n=0,1,2,\ldots,\,$
$Z_n(x^2)\equiv P_{2n}(x)$.
\end{lemma}

\noindent
\textbf{Proof.} The proof proceeds by induction. For $n=0$, by
definition, $Z_0(x^2)\equiv P_0(x)\equiv 1$  holds. For $n=1$ we have $Z_1(x^2)=x^2-1=P_{2}(x)$. Assume that
the required statement holds for $n=m-1$ and $n=m$ and show that it is true for $n=m+1$. Indeed, by \eqref{230207eq3}
and \eqref{010207eq1}, $Z_{m+1}(x^2)=(x^2-2)Z_m(x^2)-Z_{m-1}(x^2)=(x^2-2)P_{2m}(x)-P_{2m-2}(x)
=x(xP_{2m}(x)-P_{2m-1}(x))
+(xP_{2m-1}(x)-P_{2m-2}(x))-2P_{2m}(x)=xP_{2m+1}(x)+P_{2m}(x)-2P_{2m}(x)=P_{2m+2}(x)$.\epr

\begin{lemma}
\label{t_ZnExpli}
$1.$ The explicit form of the polynomial $Z_n(x)$ is
\beq
\label{010207eq2}
Z_n(x)=\sum_{i=0}^n (-1)^{i}\binom{2n-i}{i}x^{n-i}.
\eeq

$2.$ The set of roots of $Z_n(x)$ is
$\left\{4\cos^2\!{\pi k \over 2n+1}\;\Big|\;k=\1n\right\}$
and they all belong to $[0, 4[$.
\end{lemma}

\noindent
\textbf{Proof.} Due to Lemma~\ref{010207prop1}, to verify item~1, it suffices to compare  \eqref{010207eq2} with
\eqref{080207eq2}; item~2 follows from~\eqref{210507eq1}.
\epr

\smallskip
Since for $x\geq0\,$ $P_{2n}(\sqrt x)=Z_n(x)$, \eqref{200507eq2} provides the following trigonometric
representation of $Z_n(x)$:
\beq
\label{220407a}
Z_{n}(x)
={\sin((2n+1)\varphi)\over\sin\varphi},\;\;\text{where}\;\;
x=4\cos^2\varphi,\;\;\varphi\in\left]0,\tfrac{\pi}{2}\right],\;\;x\in [0,4[.
\eeq

\smallskip
\begin{lemma}
\label{270207le}
The set of the roots of the equation
\beq
\label{010707eq2}
Z_n(x)+(-1)^p=0,\quad p\in\{0,1\}
\eeq
is $\left\{4\cos^2\!{\pi k\over 2n+1+(-1)^{k+p}}\;\Big|\;k=\1n\right\}$
and they all belong to $[0, 4[$.
\end{lemma}

\medskip
\noindent
\textbf{{Proof.}} By \eqref{220407a} the roots $x_i$ of Eq.~\eqref{010707eq2} are
connected with the roots $\varphi_i$ of the equation
\beq
\label{270207a}
{\sin((2n+1)\varphi)\over\sin\varphi}+(-1)^p=0
\eeq
by $x_i=4\cos^2\!\varphi_i$.
We first solve Eq.~\eqref{270207a} and then return to \eqref{010707eq2}.
Observe that if
\beq
\label{230407b}
(2n+1)\varphi=\pi k-(-1)^{k+p}\varphi,
\eeq
where $k$ is integer, then the equality $\sin((2n + 1)\varphi)=(-1)^{p + 1}\sin\varphi$ to which
Eq.~\eqref{270207a} reduces when $\varphi\in\left]0,{\pi\over 2}\right]$ is satisfied.
By \eqref{230407b},
\beq
\label{e_roots1}
\varphi={\pi k\over 2n+1+(-1)^{k+p}}.
\eeq

Taking \eqref{e_roots1} with $k=\1n$ provides $n$ distinct roots of Eq.~\eqref{270207a}; all of them belong to
$]0,{\pi\over 2}]$. Due to \eqref{220407a}, the set of the corresponding distinct roots of Eq.~\eqref{010707eq2} in
$[0,4[$ is $\{4\cos^2{\pi k\over 2n+1+(-1)^{k+p}}\;\big|\;k=\1n\}$. Since the degree of \eqref{010707eq2} is $n$,
this equation has no other roots.
\epr

\section{Essentially cyclic digraphs with ring structure}
\label{s_main}

\x
By a digraph with ring structure we mean a digraph that contains a Hamiltonian cycle and whose remaining arcs
belong to the inverse\x Hamiltonian cycle.

More specifically, 
let $\G_n^1=(V_n,E_n^1)$ and $\G_n^2=(V_n,E_n^2)$ be the digraphs with $V_n=\{\1n\}$,
$E_n^1=\{(1,n),\,(n,n-1)\cdc(2,1)\}$, and $E_n^2=E_n^1\cup\{(1,2),\,(2,3),\dots,(n-1,n),\,(n,1)\}$. We say that
$\G_n=(V_n,E)$
is a \emph{digraph with ring structure\/} if
it is isomorphic to some $\widetilde\G_n=(V_n,\widetilde E)$ with $E_n^1\subseteq\widetilde E\subseteq E_n^2$.

\begin{figure}[htb]
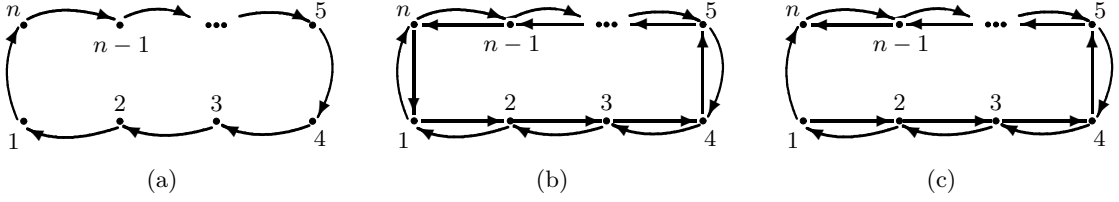

\setlength{\unitlength}{0.8mm}
\vspace{0mm}
\hspace{-1.3em}
\input g1g2g3a.lp
\begin{center}
\vspace{-2.0em}
\caption{The digraphs: (a) $\G_n^1$; (b) $\G_n^2$; (c) $\G'_n$.\label{f_g1g2g3}}
\vspace{-1.0em}
\end{center}
\end{figure}

Digraphs $\G_n^1$ and $\G_n^2$ are shown in Fig.~\ref{f_g1g2g3}(a) and \ref{f_g1g2g3}(b), respectively.

In this section, we answer the question in the title of the paper and find the Laplacian spectra of some digraphs with ring structure.
Certain \emph{weighted\/} digraphs with ring structure are considered in Section~\ref{s_weightedv EC}.

\subsection{The digraphs $\G_n^1$ with $n$ arcs and $\G_n^2$ with $2n$ arcs}

\begin{theorem}
\label{t_G1G2}
{\rm 1.} $\G_n^1$ is essentially cyclic$;$ its Laplacian spectrum is
$\left\{2\sin^2\!\frac{\pi k}{n}+\im\sin\frac{2\pi k}{n}\;\Big|\;k=\1n\right\}.$\\
{\rm 2.} $\G_n^2$ is not essentially cyclic$;$ its Laplacian spectrum is
$\left\{4\sin^2\!\frac{\pi k}{2n}\;\Big|\;k=\1n\right\}.$
\end{theorem}

\noindent \textbf{Proof.}
1. The Laplacian matrix of $\G_n^1$ is $L_n^1=I-Q$, where $Q=(q_{ij})\in\R^{n\times n}$ is the circulant permutation matrix
with the entries $q_{ij}=1$ whenever $i-j\in\{1,\,1-n\}$ and $q_{ij}=0$ otherwise. The characteristic polynomial of $Q$
is $\Delta_Q(\la)=\la^n-1$ and its spectrum is $\left\{e^{-2\pi k\im/n}\;\big|\;k=\1n\right\}$. Therefore the eigenvalues
of $L_n^1=I-Q$ are
\[
\la_k=1-e^{-2\pi k \im/n}
=1-\cos\frac{2\pi k}{n}+\im\sin\frac{2\pi k}{n}
=2\sin^2\!\frac{\pi k}{n}+\im\sin\frac{2\pi k}{n},\quad k=\1n
\]
(cf.\ \S\,{2.1} in \cite{CvetkovicDoob80}, Section 4.8.3 in \cite{MarcusMinc64}, and Section~4 in~\cite{AgaChe05LAA}).

2. The Laplacian matrix of $\G_n^2$ is symmetric and coincides with that of the undirected $n$-cycle. The spectrum
of this matrix was found in \cite{Fiedler72} and \cite{AndersonMorley} and coincides with the expression
given in Theorem~\ref{t_G1G2}. On the derivation of this expression, 
see Section~\ref{s_concl}.
\epr

The following representation of the spectrum of $\G_n^1$ is a consequence of Theorem~\ref{t_G1G2}.
\medskip\smallskip

\noindent
\textbf{Corollary of Theorem~\ref{t_G1G2}.} {\em The Laplacian spectrum of\/ $\G_n^1$ is\/
$\left\{2\sin\frac{\pi k}{n}\,e^{\pi\im\left(\frac{1}{2}-\frac{k}{n}\right)}\;\Big|\;k=\1n\right\}$}.

\subsection{The digraphs $\G'_n$ with $2n-1$ arcs}

Consider the digraph that differs from $\G_n^2$ by one arc. Let $\G'_n=(V_n,E'_n)$, where
$E'_n=E_n^2\smallsetminus\{(n,1)\}$, see Fig.~\ref{f_g1g2g3}(c).
The Laplacian matrix of $\G'_n$,
\[
L'_n=\left(%
\begin{array}{rrrrrrr}
 2  &-1     & 0     &\cdots & 0     &-1    \\
-1  & 2     &-1     &       &       & 0    \\
    &\ddots &\ddots &\ddots &       &\vdots\\
    &       &\ddots &\ddots &\ddots & 0    \\
    &       &       &-1     & 2     &-1    \\
\0L &       &       &       &-1     & 1    \\
\end{array}%
\right),
\]
differs from the tridiagonal matrix $M_n$ (see~\eqref{280609eq3}) 
by the non-zero $(1,n)$ entry.
Theorem~\ref{280207th1} below states that the Laplacian characteristic polynomial of $\G'_n$ can be expressed via the
polynomial $Z_n$ introduced in Section~\ref{s_intr} and that $\G'_n$, as well as $\G_n^2$, is not essentially cyclic.

\begin{theorem}
\label{280207th1}
Let $L'_n$ be the Laplacian matrix of 
digraph\/ $\G'_n$ whose arcs constitute the Hamiltonian cycle
$(1,n),\,(n,n-1),\dots,(2,1)$ and the path $(1,2),\,(2,3),\dots,(n-1,n),$ and there are no other arcs. Then$:$\\
\indent
{\rm 1.} The characteristic polynomial of $L'_n$ is $\Delta_{L'_n}(\la)=Z_n(\la)-(-1)^n.$ 

{\rm 2.} $\G'_n$ is not essentially cyclic and its Laplacian spectrum is
$\{4\cos^2{\pi k\over 2n+1-(-1)^{k+n}},\;k=\1n\}$.
\end{theorem}

\noindent \textbf{Proof.} 1. Expanding $\det(\la I-L'_n)$ along the first row and making use of
\eqref{010207eq1} provide
\[
\Delta_{L'_n}(\la)=\det(\la I-L'_n)=(\la-2)Z_{n-1}(\la)-Z_{n-2}(\la)-(-1)^{n}=Z_{n}(\la)-(-1)^{n}.
\]
2. By Lemma~\ref{270207le}, the roots of $\Delta_{L'_n}(\la)$ are $4\cos^2\!{\pi k\over 2n+1-(-1)^{k+n}}$,
$k=\1n$, so $\G'_n$ is not essentially cyclic.
\epr

\subsection{The digraphs $\G''_n$ with $2n-2$ arcs}

Now consider the digraphs with ring structure consisting of a Hamiltonian cycle supplemented by the
inverse Hamiltonian cycle in which \emph{two\/} arbitrary arcs are lacking. A digraph $\G''_n$ of this kind
results from $\G'_n$ by removing some $(i,i+1)$ arc, $1\le i<n$. Therefore the Laplacian matrix $L''_n$ of \/$\G''_n$ is
obtained from $L'_n$ by changing two elements in the $i$th row:
\beq
\label{L''}
L''_n=\left(%
\begin{array}{rrrrrrrrrr}
 2   &-1     &0      &\cdots &\cdots &\cdots &\cdots &0      &-1    \\
-1   & 2     &-1     &       &       &       &       &       &0     \\
     &\ddots &\ddots &\ddots &       &       &       &       &\vdots\\
     &       &\ddots &\ddots &\ddots &       &       &       &\vdots\\
     &       &       &-1     & 1     &0      &       &       &\vdots\\
     &       &       &       &\ddots &\ddots &\ddots &       &\vdots\\
     &       &       &       &       &\ddots &\ddots &\ddots & 0    \\
     &       &       &       &       &       &-1     & 2     &-1    \\
\0L  &       &       &       &       &       &       &-1     & 1    \\
   \end{array}%
\right)
\begin{array}{l}
1 \\2 \\\vdots\\\vdots \\i\\\vdots\\\vdots\\n-1\\n \\
\end{array}
\eeq

It turns out that answering the question of whether $\G''_n$ is essentially cyclic 
reduces to studying the products of Chebyshev polynomials of the second kind. To do this, we need the following notation.

Let
\beq
\label{e_x1x2}
x^{(m)}_1=4\cos^2{\pi m\over 2 m+1 }\quad\text{and}\quad x^{(m)}_2=4\cos^2{\pi(m-1)\over 2m+1}\;\;\text{(provided that\, $m>1$)}
\eeq
be the smallest and the second smallest roots of $Z_m(x)$, respectively (see Lemma~\ref{t_ZnExpli});
\beq
\label{e_u1u2}
u^{(m)}_1=4\cos^2{\pi m\over 2(m+1)}\quad\text{and}\quad u^{(m)}_2=4\cos^2{\pi(m-1)\over 2m}\;\;\text{(provided that\, $m>1$)}
\eeq
are the smallest and the second smallest roots of $Z_m(x)+(-1)^m$, respectively (see Lemma~\ref{270207le}).

\begin{figure}[htb]
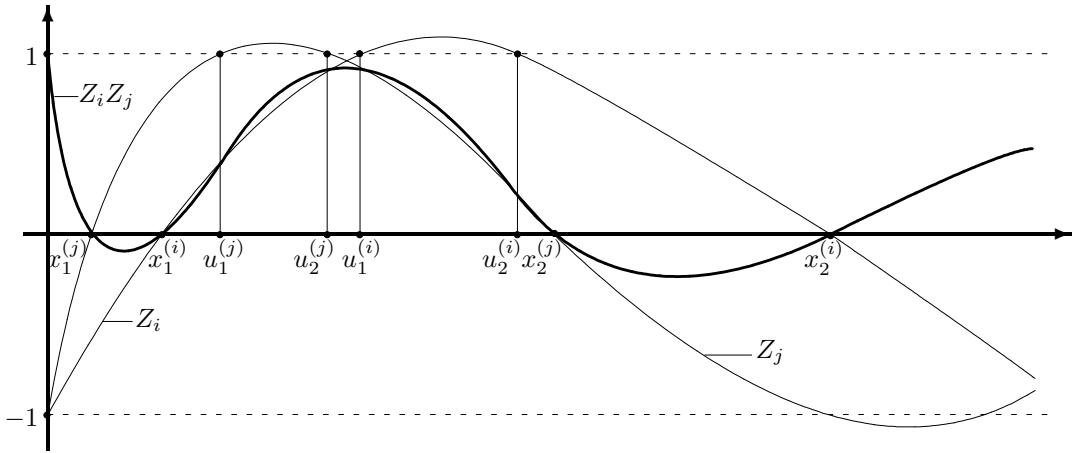

\setlength{\unitlength}{0.8mm}
\vspace{0mm}
\input raf1206.lp
\begin{center}
\vspace{-6.0em}
\caption{The shapes of the polynomials $Z_i(x)$, $Z_j(x)$, and $Z_i(x)Z_j(x)$ when $i$ and $j$ are odd and $0<i<j-1$. 
\label{f_oddcurve}}
\end{center}
\end{figure}

The following lemma establishes several inequalities involving the roots \eqref{e_x1x2}--\eqref{e_u1u2} of
$Z_m(x)$ and $Z_m(x)+(-1)^m$ as well as the polynomials that have the product form $Z_i(x)Z_j(x)$.
The shapes of $Z_i(x)$, $Z_j(x)$, and $Z_i(x)Z_j(x)$ when $i$ and $j$ are odd 
are exemplified in Fig.~\ref{f_oddcurve}.

\begin{lemma}
\label{150307th3}
Let $\:0<i<j-1$.

{\rm 1.} If $\;u_1^{(i)}< x_2^{(j)},\,$ then $\;{i\over j-1}>{2\over 3}$.

{\rm 2.} $u_1^{(i)}> u_2^{(j)}$.

{\rm 3.} The inequality
\beq
\label{150307b}
|Z_i(x)Z_j(x)|<1
\eeq
holds for all $x\in\,]0,\max(x^{(i)}_1,u_2^{(j)})]$.

{\rm 4.} If $u_1^{(i)}< x_2^{(j)},$ then \eqref{150307b} is satisfied for all $x\in [u_1^{(i)},x_2^{(j)}]$.
\end{lemma}

Lemma~\ref{150307th3} is the main technical tool employed in the proofs of the subsequent theorems of this section.

\medskip
\noindent
\textbf{Proof.} {\rm 1.} Since $\cos^2 t$ strictly decreases on $[0,{\pi\over 2}]$, \,
$u_1^{(i)}<x_2^{(j)}$ implies that ${i \over 2(i+1)}> {j-1 \over 2j+1}$. As a consequence, one
has item~1 of Lemma~\ref{150307th3}.

Item 2 follows from the definitions of $u_1^{(i)}$ and $u_2^{(j)}$ and the inequality $0<i<j-1$.

{\rm 3.} Let $x_1^{(i)}\geq u_2^{(j)}$. Let us show that  \eqref{150307b} holds for all   $x\in\,]0,x_1^{(i)}]$.
If $x\in\,]0,x_1^{(i)}[$, then $\varphi\in\,]{\pi i\over 2i+1},{\pi\over 2}[$, where $\varphi=\arccos\frac{x}{2}$.
Let $\varphi={\pi(i+\delta)\over 2i+1}$ for $0<\delta<1/2$. We get
\beq
\label{280307eq1}
|Z_i(x)Z_j(x)|
\leq\Big|{\sin(2i+1)\varphi\over\sin^2\varphi}\Big|
=  \left|{\sin\pi\delta    \over\sin^2{\pi(i+\delta)\over 2i+1}}\right|.
\eeq
Consider the derivative of\, $\sin\delta\pi\big\slash\sin^2{\pi(i+\delta)\over 2i+1}$:
\begin{eqnarray}
\nonumber  \left({\sin\pi\delta\over\sin^2{\pi(i+\delta)\over 2i+1}}\right)'_{\delta}
           &=&{\pi\cos(\pi\delta)\sin{\pi(i+\delta)\over 2i+1}-{2\pi\over 2i+1}\sin(\pi\delta)\cos{\pi(i+\delta)\over
2i+1} \over \sin^3{\pi(i+\delta) \over 2i+1}}\\
\nonumber  &=&{\big({\pi(2i-1)\over 2i+1}+{2\pi\over 2i+1}\big)\cos(\delta\pi)\sin{\pi(i+\delta)\over 2i+1}-{2\pi\over
2i+1}\sin(\delta\pi)\cos{\pi(i+\delta)\over 2i+1}\over\sin^3{\pi(i+\delta)\over 2i+1}}\\
\label{240407a}
           &=&{{\pi(2i-1)\over 2i+1}\cos(\pi\delta)\sin{\pi(i+\delta)\over 2i+1}+{2\pi\over 2i+1}
\sin{\pi i(1-2\delta)\over 2i+1}\over\sin^3{\pi(i+\delta)\over 2i+1}}.
\end{eqnarray}

Since for $0\leq\delta<1/2$ the trigonometric functions in \eqref{240407a} are positive,
so is the left-hand side, therefore        ${\sin\pi\delta \big / \sin^2{\pi(i+\delta) \over 2i+1}}$ increases in $\delta$.
For $\delta=0$ and $\delta=1/2$, ${\sin\pi\delta \big / \sin^2{\pi(i+\delta) \over 2i+1}}$ equals $0$
and $1$, respectively. Hence
${\sin\pi\delta\big/\sin^2{\pi(i+\delta)\over 2i+1}}<1$ holds for all $\delta\in[0,{1\over2}[$ and, by \eqref{280307eq1},
$|Z_i(x)Z_j(x)|<1$ is satisfied for all $x \in\,]0,x_1^{(i)}]$.

Let us show that $u_2^{(j)}>x_1^{(i)}$ implies $|Z_i(x) Z_j(x)|<1$  as well. Consider the parametrization
$\varphi=\varphi(\delta)={\pi \over 2}{j-\delta \over j}$, from which
$x=x(\delta)=4\cos^2{\pi \over 2}{j-\delta \over j}$. Then we obtain $x(0)=0,\, x\big(\!\frac{j}{j+1}\big)=u_1^{(j)},$ and
$x(1)=u_2^{(j)}$. For $x\in\,]0,u_2^{(j)}]$  (i.\,e., for  $\delta\in\,]0,1]$) we have
\[
Z_i(x)Z_j(x)=
{\sin\!\big(\!{\pi\over 2}(2i+1){j-\delta\over j}\big)\sin  \!\big(\!{\pi\over 2}(2j+1){j-\delta\over j}\big)\over
                                                      \sin^2\!\big(\!{\pi\over 2}      {j-\delta\over j}\big)},
\]
consequently,
\[
|Z_i(x)Z_j(x)|\leq
 {\big|\sin  \!\big(\!{\pi\over 2}(2i+1)\big(1-{\delta\over j}\big)\!\big)\!\big|\over
       \sin^2\!\big(\!{\pi\over 2}\big(1-{\delta\over j}\big)\!\big)}
={\big|\cos  \!\big(\!{\delta\pi\over 2} {2i+1  \over j}\big)\!\big|\over
       \cos^2\!\big(\!{\delta\pi\over 2} {   1  \over j}\big)}.
\]

If $\cos\!\big(\!{\delta\pi\over 2}{2i+1\over j}\big)>0$, then ${\delta\pi\over 2}{2i+1\over j}<{\pi\over 2}$, so
\[
|Z_i(x)Z_j(x)|\leq
{ \cos  \!\big(\!{\delta\pi\over 2}{2i+1\over j}\big)\over\cos^2\!\big(\!{\delta\pi\over 2}{1\over j}\big)}
<{\cos  \!\big(\!{\delta\pi\over 2}{2i  \over j}\big)\over\cos^2\!\big(\!{\delta\pi\over 2}{1\over j}\big)}
={\cos^2\!\big(\!{\delta\pi\over 2}{ i  \over j}\big)-    \sin^2\!\big(\!{\delta\pi\over 2}{i\over j}\big)\over
  \cos^2\!\big(\!{\delta\pi\over 2}{1   \over j}\big)}
<1.\]

If $\cos\!\big(\!{\delta\pi \over 2}{2i+1 \over j}\big)\leq 0$, then ${\pi\over 2}\leq{\delta\pi\over 2}{2i+1\over j}<\pi$,
and since $i+1<j$,                                           we have ${\pi\over 2}\leq{      \pi\over 2}{2i+1\over j}<\pi$. Therefore
\begin{eqnarray*}
      {\big|\!\cos  \!\big(\!{\delta\pi\over 2}{    2i+1\over j}\big)\!\big|\over
              \cos^2\!\big(\!{\delta\pi\over 2}{       1\over j}\big)}
&\leq&{\big|\!\cos  \!\big(\!{      \pi\over 2}{    2i+1\over j}\big)\!\big|\over
              \cos^2\!\big(  {      \pi\over 2}{       1\over j}\big)}
=            {\cos  \!\big(  \pi-  {\pi\over 2}{    2i+1\over j}\big)\over
              \cos^2\!\big(  {      \pi\over 2}{       1\over j}\big)}
=            {\cos  \!\big(  {      \pi\over 2}{2(j-i)-1\over j}\big)\over
              \cos^2\!\big(  {      \pi\over 2}{       1\over j}\big)}\\
&\leq&       {\cos  \!\big(\!{      \pi\over 2}{       2\over j}\big)\over
              \cos^2\!\big(\!{      \pi\over 2}{       1\over j}\big)}
=            {\cos^2\!\big(\!{      \pi\over 2}{       1\over j}\big)-
              \sin^2\!\big(\!{      \pi\over 2}{       1\over j}\big)\over
              \cos^2\!\big(\!{      \pi\over 2}{       1\over j}\big)}
<1.
\end{eqnarray*}

4. By item~1 of Lemma~\ref{150307th3}, if $u_1^{(i)}< x_2^{(j)}$, then $i> {2j-2 \over 3}$. On the other hand,
$i+1<j$, and due to these two inequalities,
\vspace{-.5em}
\beq
\label{ijge8}
i+j\geq 8
\eeq
\vspace{-.5em}
holds. The segment $\big[{\pi(j-1) \over 2j+1},{\pi i\over 2(i+1)}\big]$ on the $\varphi$-axis corresponds to
$x\in[u_1^{(i)},x_2^{(j)}]$. Let $\varphi({\xi})={{\pi(j-1+\xi)\over 2j+1}}$, where $0\leq\xi\leq {3i-2j+2
\over 2i+2}$. Then $\varphi(0)={{\pi(j-1) \over 2j+1}}$ and $\varphi\big({3i-2j+2 \over 2i+2}\big)={\pi i\over 2(i+1)}$.

For all $x\in [u_1^{(i)}, x_2^{(j)}]$, the following inequality holds:
\begin{eqnarray}
\nonumber
|Z_i(x)Z_j(x)|
&\leq&{    |\sin        (  (2j+1)\varphi( \xi)                )      |\over\sin^2\varphi(  \xi)}
 =    {\big|\sin  \!\big(\!(2j+1){\pi(j-1+\xi) \over 2j+1}\big)\!\big|\over\sin^2 {\pi(j-1+\xi)\over 2j+1}}\\
&=   &{    |\sin        (         \pi(j-1+\xi)                )      |\over
            \sin^2\!\big({\pi\over 2}-{\pi({3 \over 2}-\xi    )       \over 2j+1}\big)}
 =    {    |\sin        ( \pi             \xi                 )      |\over
            \cos^2       {\pi({3\over 2}- \xi)\over 2j+1}}.
\label{060507a}
\end{eqnarray}
Since  $i<j-1$, we get
$\xi\leq{3i-2j+2\over 2i+2}
={1\over 2}\big({5i+5-2(i+j)-3\over 2i+2}        \big)
={1\over 2}\big( 5  -{2(i+j)+3\over  i+1}        \big)\leq
 {1\over 2}\big( 5  -{2(i+j)+3\over {i+j\over 2}}\big)
={i+j-6\over 2(i+j)}<{1\over 2}$ and \eqref{060507a} results in 
\[
|Z_i(x)Z_j(x)|
\leq{    |\sin(\pi\xi)|                                         \over\cos^2{ \pi({3\over 2}-\xi)\over 2j+1}}
\leq{\big|\sin\!\big(  \pi{i+j-6\over 2(i+j)}\big)\!       \big|\over\cos^2{     {3\over 2}\pi  \over  i+j}}
<   {\big|\cos\!\big(  {3\pi    \over   i+j }\big)\!       \big|\over\cos^2{3\pi   \over 2(i+j)}}
=   {\big|\cos^2{3\pi\over 2(i+j)}-\sin^2{3\pi\over 2(i+j)}\big|\over\cos^2{3\pi   \over 2(i+j)}}
<1.
\]
The last inequality follows from \eqref{ijge8}. The lemma is proved.
\epr
\medskip\smallskip

In all subsequent statements, the copies of multiple roots of a polynomial are considered as distinct roots so that
each polynomial of degree $n$ has $n$ distinct roots.

\begin{lemma}
\label{t_ZZle1}
Let $\:0<i<j-1$. Let $x_1,\,x_2,$ and $x_3$ be the three smallest
roots of the polynomial $f(x)=Z_i(x)Z_j(x)$ and $x_1<x_2\le x_3$. Then $|f(x)|<1$ for all $x\in\,]0,x_3].$
\end{lemma}

\noindent
\textbf{Proof.}
Consider four cases.
(a) $x^{(i)}_1<   u_2^{(j)}$ and $u_1^{(i)}<   x_2^{(j)}$.
This case is illustrated by Fig.~\ref{f_oddcurve}. By item~3 of Lemma~\ref{150307th3}, $|f(x)|<1$ holds for all
$x\in\,]0,u_2^{(j)}]$; by item~4, this inequality is also true on $[u_1^{(i)},x_2^{(j)}]$. By item~2, $u_2^{(j)}<u_1^{(i)}$.
On $]u_2^{(j)},u_1^{(i)}[$, we also have $|f(x)|<1$, because on this interval $|Z_i(x)|<1$ and $|Z_j(x)|<1$. Thus,
$|f(x)|<1$ on $]0,x_2^{(j)}]$. Since $x_2^{(j)}=x_3$, the desired statement follows.

(b) $x^{(i)}_1<   u_2^{(j)}$ and $u_1^{(i)}\ge x_2^{(j)}$.
In this case, using item~3 of Lemma~\ref{150307th3}, we similarly obtain $|f(x)|<1$ on $]0,u_1^{(i)}]$.
Since $0<x_3=x_2^{(j)}\le u_1^{(i)}$, $|f(x)|<1$ is true on $]0,x_3].$

(c) $x^{(i)}_1\ge u_2^{(j)}$ and $u_1^{(i)}<   x_2^{(j)}$.
By items~3 and~4 of Lemma~\ref{150307th3}, $|f(x)|<1$ holds on $]0,x^{(i)}_1]$ and $[u_1^{(i)},x_2^{(j)}]$. In addition,
$|f(x)|<1$ on $]x^{(i)}_1,u_1^{(i)}[$, because on this interval, $|Z_i(x)|<1$ and $|Z_j(x)|<1$. Hence $|f(x)|<1$ on
$]0,x_2^{(j)}]$, where $x_2^{(j)}=x_3$.

(d) $x^{(i)}_1\ge u_2^{(j)}$ and $u_1^{(i)}\ge x_2^{(j)}$.
In this case, $|f(x)|<1$ is guaranteed on $]0,x^{(i)}_1]$. If $x^{(i)}_1\ge x^{(j)}_2$, then $x_3=x^{(i)}_1$ and the
desired statement follows. In the opposite case, $x_3=x^{(j)}_2$ and $|f(x)|<1$ holds on $]x^{(i)}_1,x^{(j)}_2]$, because
on this interval $|Z_i(x)|<1$ and $|Z_j(x)|<1$. Therefore $|f(x)|<1$ on $]0,x_3]$, as needed.
\epr

The next lemma provides a means of using Lemma~\ref{t_ZZle1} in the subsequent proofs.

\begin{lemma}
\label{t_poly3}
Suppose that $g(x)$ is a polynomial with real coefficients and $x_1\le x_2\le x_3$ are some of its distinct real roots.
Suppose that $|g(x)|<1$ for all $x$ such that $x_1\le x\le x_3$. Then each of the polynomials $g(x)-1$ and $g(x)+1$ has
at least a pair of non-real roots.
\end{lemma}

\noindent
\textbf{Proof.} Under the assumptions of Lemma~\ref{t_poly3},
neither $g(x)+1$ nor $g(x)-1$ has any real root $x'$ such that $x_1\le x'\le x_3$. On the other hand,
the segment (possibly, degenerating into a point) $[x_1,x_3]$ contains at least two roots of the derivative $g'(x)$.
Consequently, each of the polynomials $g(x)+1$ and $g(x)-1$ has no root non-strictly between two roots of its derivative.
Hence it has at least two non-real roots.
\epr

The following theorem determines which digraphs of type $\G''_n$ are essentially cyclic. Its proof is based on
Lemmas~\ref{t_ZZle1} and~\ref{t_poly3} and the subsequent Lemma~\ref{010707le1}.

\begin{theorem}
\label{280207th2}
Let $L''_n$ be the Laplacian matrix of the digraph\/ $\G''_n$ whose arcs form the Hamiltonian cycle
$(1,n),\,(n,n-1),\dots,(2,1),$ the path $(1,2),\,(2,3),\dots,(i-1,i),$ and the path $(i+1,i+2),\dots,(n-1,n),$ where $1\le i<n$.
Then$:$

{\rm 1}. The characteristic polynomial of $L''_n$ is $\Delta_{L''_n}(\la)=Z_i(\la)Z_{n-i}(\la)-(-1)^{n}.$

{\rm 2}. If $n$ is even$,$ then
\/$\G''_n$ is essentially cyclic for all $i\in\{1\cdc n-1\}$ except for \,$i={n\over 2};$
in the latter case the eigenvalues of $L''_n$ are\/
$4\cos^2{\pi k\over n  } $ and\/
$4\cos^2{\pi k\over n+2},\:k=1\cdc{n\over 2}$.

{\rm 3}. If $n$ is odd$,$ then \/$\G''_n$ is essentially cyclic for all $i\in\{1\cdc n-1\}$ except for \,$i={n-1\over 2}$ and $i={n+1\over 2};$
in the latter cases the eigenvalues of $L''_n$ are\/ $4\cos^2\!{\pi k\over n+1},\:k=1\cdc n.$
\end{theorem}

\begin{figure}[htb]
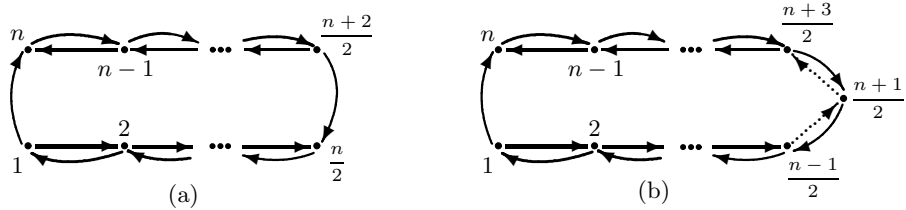

\setlength{\unitlength}{0.8mm}
\addvspace{0mm}\vspace{0mm}
\hspace{2em}
\input temp1.lp
\begin{center}
\vspace{-1.8em}
\caption{Digraphs $\G''_n$ that are not essentially cyclic:
(a) $n$ is even; $i=n/2$;
(b) $n$ is  odd; either $i=(n-1)/2$ or $i=(n+1)/2$, i.\,e., exactly one of the two dotted vectors is an arc.\label{f_uness}}
\vspace{-0.6em}
\end{center}
\end{figure}

The only digraphs of type $\G''_n$ 
that are not essentially cyclic have their two vertices of indegree~1 (as well as the two vertices of outdegree~1) at a maximum possible distance. In other words, they can be obtained from $\G_n^2$ by the removal of two \emph{most distant\/} arcs from the same Hamiltonian cycle. These digraphs are shown in Fig.~\ref{f_uness}. In Fig.~\ref{f_uness}(b), exactly one of the two
dotted vectors must be an arc; the two resulting digraphs are obviously isomorphic.

Fig.~\ref{f_shapes} shows the shape of the polynomials $Z^2_i(\la)$ and $Z_i(\la)Z_{i+1}(\la)$, which differ by $1$ from the characteristic polynomials $\Delta_{L''_n}(\la)$ of the digraphs $\G''_n$ that are not essentially cyclic due to Theorem~\ref{280207th2}.

\begin{figure}[htb]
\vspace*{0mm}
\begin{center}
\setlength{\unitlength}{0.8mm}
\hspace{0mm}
\input raf2804.lp
\vspace{-1em}
\caption{(a) The shape of $Z^2_i(\la)$; (b) the shape of $Z_i(\la)Z_{i+1}(\la)$.\label{f_shapes}}
\vspace{-1.0em}
\end{center}
\end{figure}

\medskip
\noindent
\textbf{Proof of Theorem~\ref{280207th2}.} 1. Expanding $\Delta_{  L''_n}(\la)=\det(\la I-L''_n)$ along the
first row and using identity \eqref{010207eq1} and the fact that for any square matrices $P$ and $S$,
$\det\!\left(\begin{array}{c|c} P&0\\\hline R&S\\\end{array}\right)\!=\det P\det S$, one obtains
\[
\Delta_{  L''_n}(\la)
= (\la-2)Z_{i-1}  (\la)\,Z_{n-i}(\la)-Z_{i-2}(\la)\,Z_{n-i}(\la)-(-1)^{n}
=        Z_i      (\la)  Z_{n-i}(\la)                           -(-1)^{n},
\]
which proves item~1 of Theorem~\ref{280207th2}. To prove items~2 and~3, we need the following lemma.

\begin{lemma}
\label{010707le1}
{\rm 1.}
If\/ $i+j$ is even$,$ then the equation $Z_i(x)Z_j(x)-1=0$ has only real roots if and
only if $i=j$. In the latter case$,$ the roots are\/ $4\cos^2{\pi k\over 2j},$ $4\cos^2{\pi k\over 2j+2},\,k=1\cdc j$.

{\rm 2.} If $i+j$ is odd and $i<j,$ then the equation
$Z_i(x)Z_j(x)+1=0$  has only real roots if and only if $i=j-1$. In the latter case$,$
\beq
\label{250407a}
Z_i(x)Z_j(x)+1=P^2_{i+j}\big(\sqrt{x}\big)\;\;\text{for all}\;\;\:x\ge0
\eeq
holds and the roots are\/ $4\cos^2\!{\pi k\over 2j},\:k=1\cdc i+j.$
\end{lemma}

\noindent
\textbf{Proof of Lemma \ref{010707le1}.} We first prove that under the requirements of Lemma~\ref{010707le1}, all
the roots are real. After that we show that otherwise there are at least two non-real roots.

1. If  $i=j$, then $Z_i(x) Z_j(x)-1=(Z_i(x)-1)(Z_i(x)+1)$, thus, Lemma~\ref{270207le} implies that all the roots
are real, belong to $[0,\,4[$, and can be expressed as\, $4\cos^2{\pi k\over 2j}$, $4\cos^2{\pi k\over 2j+2},\;k=1\cdc j$.
This can also be obtained using Lemma~\ref{010207prop1}, Eq.~\eqref{210507eq1}, and Catalan's identity for Chebyshev polynomials (see \cite[Eq.~$(1.1')$]{Udrea95}) $(P_n(x)-1)(P_n(x)+1)=P_{n-1}(x)\,P_{n+1}(x)$.

2. Let  $j=i+1$. Then for $x\in [0,4[$ and $\varphi=\arccos{\sqrt x\over 2}$, using \eqref{220407a}, \eqref{080207eq1},
and \eqref{e_P2m1}, we have
\begin{eqnarray*}
Z_i(x)Z_j(x)+1
&=& {\sin((2i+1)\varphi) \sin((2i+3)\varphi) \over \sin^2\varphi}+1
={\cos(2\varphi)-\cos((4i+4)\varphi)\over 2\sin^2\varphi}+1\\
\nonumber
&=&{\sin^2((2i+2)\varphi)\over\sin^2\varphi}
=P^2_{i+j}\big(\sqrt{x}\big)
=x\prod_{k=1}^{j-1}\Big(x-4\cos^2\!{\pi k\over 2j}\Big)^{\!2}
= \prod_{k=1}^{i+j}\Big(x-4\cos^2\!{\pi k\over 2j}\Big).
\end{eqnarray*}

This can also be obtained using the identity mentioned in the proof of item~1.
Thus, the roots of $Z_i(x)Z_j(x)+1$ are $4\cos^2\!{\pi k\over 2j},\:k=1\cdc i+j$;
they are real and belong to $[0,\,4[$.

Let us prove that in the remaining cases, each of the equations under consideration has at least a pair of non-real roots.

Let $0<i<j-1$.
By Lemma~\ref{t_ZZle1}, $|Z_i(x)Z_j(x)|<1$ for all $x\in\,]0,x_3]$, where
$x_1,\,x_2,$ and $x_3$ ($x_1<x_2\le x_3$) are the three smallest roots of $Z_i(x)Z_j(x)$.
Hence by Lemma~\ref{t_poly3} each of the polynomials $Z_i(x)Z_j(x)+1$ and $Z_i(x)Z_j(x)-1$ has at least a pair of
non-real roots.
The lemma is proved.
\epr

In view of item~1 of Theorem~\ref{280207th2}, items~2 and~3 follow from items 1 and 2 of Lemma \ref{010707le1}, respectively.
\epr

\subsection{The digraphs $\G_n$ with $m$ ($n<m<2n-2$) arcs}

Let us summarize the above results. According to Theorems~\ref{t_G1G2} to \ref{280207th2}: $\G_n^1$ is essentially cyclic; the digraphs $\G''_n$ are essentially cyclic except for the cases specified in Theorem~\ref{280207th2}; $\G'_n$ and $\G_n^2$ are not essentially cyclic. The following theorem answers the question of essential cyclicity for the remaining digraphs with ring structure, in which the Hamiltonian counter-cycle lacks more than two arcs. According to this theorem, all such digraphs are essentially cyclic.

\begin{theorem}
\label{t_gener}
Let\/ $\G_n$ be a digraph on $n>3$ vertices constituted by the Hamiltonian cycle $\{(1,n),\,(n,n\!-\!1)\cdc(2,1)\}$
                                                          and the opposite cycle $\{(1,2),\,(2,3),\dots,(n-1,n),\,(n,1)\}$
in which $i$ $(2<i<n)$ arbitrary arcs are missing. Then\/$:$ 

{\rm 1.} The characteristic polynomial of the Laplacian matrix\/ $L_n$ of\/ $\G_n$ is
\vspace{-.4em}
\beq
\label{e_DL}
\Delta_{L_n}(\la)=\prod_{k=1}^K Z_{i_k}(\la)-(-1)^{n},
\eeq

\vspace{-.4em}\noindent
where 
$i_1\cdc i_K$ are the path lengths in the decomposition of the cycle $\{(1,n),\,(n,n\!-\!1)\cdc(2,1)\}$ into the paths linking the consecutive
vertices of indegree~$1$ in $\G_n$.

{\rm 2.} $\G_n$ is essentially cyclic.
\end{theorem}

\noindent
\textbf{Proof.}
1. The proof is quite similar to that of item~1 of Theorem~\ref{280207th2}.

2. We need the following lemma, which extends Lemma~\ref{t_ZZle1}.

\begin{lemma}
\label{t_prodZ}
Let $f(x)=\prod_{k=1}^K\!Z_{i_k}(x),$ where $K>2$ and $i_k>0,\;k=1\cdc K.$

$1.$ Let $x_1,\,x_2,$ and $x_3$ be the 
smallest\/\footnote{Here, as earlier, we distinguish the copies of every multiple root of a polynomial.}
roots of $f(x)$ and {$x_1\!\le\!x_2\!\le\! x_3$.} Then\/ $|f(x)|\!<\!1$ for all $x\in\,]0,x_3].$

$2.$ Each of the polynomials $f(x)-1$ and $f(x)+1$ has at least a pair of non-real roots.
\end{lemma}

\noindent
\textbf{Proof of Lemma~\ref{t_prodZ}.}
1. The proof proceeds by induction on $K$. We first consider the step of induction and then come back to its base. Assume that
the required statement is true for all $K<K_0$. Let us prove that it is also true for $K=K_0$.
Consider any product of the form $f(x)=\prod_{k=1}^{K_0}\!Z_{i_k}(x).$ Without loss of generality, assume that
\beq
\label{e_iK}
i_{K_0}=\min(i_1\cdc i_{K_0}).
\eeq
Then
\beq
\label{e_fx}
f(x)=Z_{i_{K_0}}(x)\,f_0(x),
\eeq
where $f_0(x)=\prod_{k=1}^{K_0-1}\!Z_{i_k}(x)$.
Let $x_1^0,\,x_2^0,$ and $x_3^0$ be the three smallest roots of $f_0(x)$ and $x_1^0\le x_2^0\le x_3^0$.
By the assumption,\footnote{The inequality $|f_0(x)|<1$ is given here in a weakened form for the subsequent use of
the induction step in the case where the strict inequality is not satisfied.}
\beq
\label{e_f01}
|f_0(x)|\le1\;\;\text{on}\;\;]0,x_3^0].
\eeq
By \eqref{e_x1x2}, the smallest root of $Z_i(x)$ decreases with the increase of~$i$. Therefore \eqref{e_iK} implies that
$x_1=x_1^0$, $x_2=x_2^0$, and $x_3=\min(x_3^0,x^{(K_0)})$, where $x^{(K_0)}$ is the smallest root of $Z_{i_{K_0}}(x)$.
Having in mind \eqref{e_f01} and that $|Z_{i_{K_0}}(x)|<1$ on $]0,x^{(K_0)}]$ (which follows from Lemma~\ref{010207prop1}
and~\eqref{080207eq1}), we have that $|f_0(x)|\le1$ and $|Z_{i_{K_0}}(x)|<1$ on $]0,x_3]$.
Hence, by \eqref{e_fx},\, $|f(x)|<1$ on $]0,x_3]$, thus, the induction step is complete.

We now turn to the base of induction. Let $K=3$ and $f(x)=\prod_{k=1}^3 Z_{i_k}(x).$ Without loss of generality, assume
that $i_3\le i_2\le i_1$. Then $f(x)=Z_{i_3}(x)f_0(x),$ where $f_0(x)=Z_{i_2}(x)\,Z_{i_1}(x)$.
Let $x_1^0,\,x_2^0$, and $x_3^0$ be the three smallest roots of $f_0(x)$ ordered as follows: $x_1^0\le x_2^0\le x_3^0$.
Consider three cases.

(a) $i_2<i_1-1$. Then, by Lemma~\ref{t_ZZle1},  $|f_0(x)|<1$ for all $x\in\,]0,x_3^0]$.
Now, applying the above induction step, one has $|f  (x)|<1$ for all $x\in\,]0,x_3  ]$, as needed.

(b) $i_2=i_1-1$. In this case, by \eqref{250407a}, $f_0(x)=Z_{i_2}(x)Z_{i_1}(x)=P^2_{i_2+i_1}\big(\sqrt{x}\big)-1$ (see
also Fig.~\ref{f_shapes}(b)) and we only have $|f_0(x)|\le1$ on $]0,x_3^0]$. However, the case of the weak inequality
$|f_0(x)|\le1$ is covered by the above induction step (see \eqref{e_f01}), thereby, this step provides $|f(x)|<1$ on $]0,x_3]$.

(c) $i_2=i_1$. In this case, by Lemma~\ref{010207prop1}, $f_0(x)=Z_{i_2}^2(x)=P_{2i_2}^2(\sqrt{x})$, and $|f_0(x)|>1$
is possible for some $x\in\,]0,x_3^0]$ (cf.\ Fig.~\ref{f_shapes}(a)). Let $x^{(i_3)}$ be the smallest root of $Z_{i_3}(x)$.
Then by \eqref{200507eq2} for every integer $k>0$ and $x\in\,]0,1]$ we have
\beq
\label{e_P2k}
P_{2k}^2(\sqrt{x})
=\frac{\sin^2\bigl((2k+1)\arccos\frac{\sqrt{x}}{2}\bigr)}{1-\frac{x}{4}}
\le\frac{1}{1-\frac{x}{4}}
=1+\frac{x}{4-x}
\le 1+\frac{x}{3}.
\eeq

Moreover, since $x^{(i_3)}\le1$,
\vspace{-.3em}
\beq
\label{e_Zi3}
|Z_{i_3}(x)|\le|Z_1(x)|=1-x\;\;\;\text{for all}\;\;\;x\in\,]0,x^{(i_3)}].
\eeq
\vspace{-.2em}
Indeed, $|Z_{i_3}(0)|=|Z_1(0)|=1$; by \eqref{010207eq2}, $|Z_{i_3}(x)|'_{x=0}=-\frac{n(n+1)}{2}\le-1=|Z_1(x)|'_{x=0}.$
Now the assumption that $|Z_{i_3}(x)|>|Z_1(x)|$ at some $x\in\,]0,x^{(i_3)}]$
implies that $|Z_{i_3}(x)|$ has an inflection on $]0,x^{(i_3)}[$, which is impossible because $Z_{i_3}(x)$ has $i_3$ real
roots and $x^{(i_3)}$ is the smallest one.

Using \eqref{e_P2k} and \eqref{e_Zi3} for every $x\in\,]0,x^{(i_3)}]$ we have 
\[
|f(x)|=|Z_{i_3}(x)\,P_{2i_2}^2(\sqrt{x})|<(1-x)(1+x)<1.
\]
Finally, $i_3\le i_2=i_1$ implies that $x_3=x^{(i_3)}$, thereby $|f(x)|<1$ for all $x\in\,]0,x_3]$.

2. Item 2 follows from item 1 and Lemma~\ref{t_poly3}. Lemma~\ref{t_prodZ} is proved.
\epr

Now item 2 of Theorem~\ref{t_gener} follows from \eqref{e_DL} and item~2 of Lemma~\ref{t_prodZ}.
\epr

\medskip
\noindent
\textbf{Corollary of Theorem~\ref{t_gener}.} {\em If for two digraphs of the type described in Theorem~$\ref{t_gener},$
the path lengths $i_1\cdc i_K$ in the decomposition of\/ $\{(1,n),\,(n,n\!-\!1)\cdc(2,1)\}$ into the paths linking the consecutive vertices of indegree~$1$ differ only by the order of the corresponding paths in the decomposition$,$ then these digraphs have the same Laplacian spectrum.}
\medskip\medskip

Thus, the answer to the question in the title of this paper is as follows.
The digraphs $\G_n$ with ring structure, which consist of two opposite Hamiltonian cycles from one of which
some arcs can be removed, are essentially cyclic, except for (up to isomorphism) three digraphs:
\begin{itemize}\vspace{-.4em}
\item $\G_n^2$, where no  arcs are removed (Fig.~\ref{f_g1g2g3}(b));\vspace{-.3em}
\item $\G'_n$,  where one arc  is  removed (Fig.~\ref{f_g1g2g3}(c)); \vspace{-.3em}
\item $\G''_n$  with two most distant\footnote{For the exact meaning of ``most distant,'' see Theorem~\ref{280207th2}.} arcs removed (Fig.~\ref{f_uness}).
\end{itemize}

A by-product of this work is the following theorem on the Chebyshev polynomials of the second kind.

\begin{theorem}
\label{t_Cheby1}
Let $h(x)=\prod_{k=1}^K\!P_{2i_k}(x)+(-1)^p,$ where $p\in\{0,1\},$ $P_{2i_k}(x)$ are the Chebyshev polynomials of the second kind scaled on $]\!-\!2,2[$ $($see \eqref{200507eq2}$),$ $K\ge1,$ and $i_k>0,\;k=1\cdc K.$
Then $h(x)$ has only real roots if and only if

{\rm (a)} $K=1;$                           the roots are $\bigl\{\pm2\cos{\pi k\over 2j+1+(-1)^{k+p}}\;\big|\;k=1\cdc j\bigr\},$ where $j=i_1$ or

{\rm (b)} $K=2,$ $i_1=i_2,$ and $p=1;$     the roots are $\bigl\{\pm2\cos{\pi k\over 2j};\;\pm2\cos{\pi k\over 2j+2}\;\big|\;k=1\cdc j\bigr\},$ where $j=i_1$ or

{\rm (c)} $K=2,$ $|i_1-i_2|=1,$ and $p=0;$ the roots are\/\footnote{In this expression, each element appears twice, which corresponds to multiplicity~2 of every root of $h(x)$ in the case~(c).}
                                                         $\bigl\{\pm2\cos{\pi k\over 2j}\,\big|\,k=1\cdc 2j-1\bigr\},$ where $j=\max(i_1,i_2)$.
\end{theorem}

\noindent
\textbf{Proof.}
By Lemma~\ref{010207prop1}, $h(x)=\prod_{k=1}^K\!Z_{i_k}(x^2)+(-1)^p$.
In the case (a), $h(x)$ has only real roots by virtue of Lemma~\ref{270207le}.

Suppose that $K=2$.

If $|i_1-i_2|>1$, then Lemmas~\ref{t_ZZle1} and~\ref{t_poly3} imply that $h(x)$ has at least a pair of non-real roots.

If $ i_1=i_2   $ and $p=1$ (the case~(b)), then by item~1 of Lemma~\ref{010707le1}, $h(x)$ has only real roots.

If $ i_1=i_2   $ and $p=0$, then $h(x)=Z_{i_1}^2(x^2)+1$ has no real roots.

If $|i_1-i_2|=1$ and $p=0$ (the case~(c)), then by item~2 of Lemma~\ref{010707le1}, $h(x)$ has only real roots.

If $|i_1-i_2|=1$ and $p=1$, then using the notation $i=\min(i_1,i_2)$, by item~2 of Lemma~\ref{010707le1} we have
$
h(x)=Z_i(x^2)Z_{i+1}(x^2)-1=P_{2i+1}^2(x)-2. 
$
Since $P_{2i+1}^2(x)$ has a maximum on $]0,1[$ and on this interval,
\[
P_{2i+1}^2(x)
=\frac{\sin^2\bigl((2i+2)\arccos\frac{x}{2}\bigr)}{1-\frac{x^2}{4}}
<\frac{1}{1-\frac{1}{4}}
=\frac{4}{3},
\]
$h(x)$ has a negative maximum, consequently, it has at least a pair of non-real roots.

Finally, if $K>2$, then by item~2 of Lemma~\ref{t_prodZ}, $h(x)$ has at least two non-real roots.

The expression for the roots in the case (a) is provided by Lemma~\ref{270207le}; in the cases (b) and (c) they are easily obtained using the identity $P_{n-1}(x)\,P_{n+1}(x)+1=P_n^2(x)$ (see \cite[Eq.~$(1.1')$]{Udrea95}) or Lemma~\ref{010707le1}.
\epr

\section{On the essential cyclicity of weighted digraphs} 
\label{s_weightedv EC}

In this section, we study the essential cyclicity of simple weighted digraphs with ring structure.

Recall that the Laplacian matrix of a weighted digraph $\G$ with strictly positive arc weights is the matrix
$L=L(\G)=(\l_{ij})$ in which, for $j\ne i$, $\l_{ij}$ equals minus the weight of arc $(i,j)$ in $\G$ and $\l_{ij}=0$
if $\G$ has no $(i,j)$ arc, the diagonal entries of $L(\G)$ being such that the row sums are zero. Some spectral
properties of the Laplacian matrices of weighted digraphs were studied in \cite{AgaChe00,AgaChe01,AgaChe05LAA,CheAga02,CauVee06Comb}, papers cited therein, and, with a different definition of the Laplacian matrix, in~\cite{Chung05}.

\begin{figure}[htb]
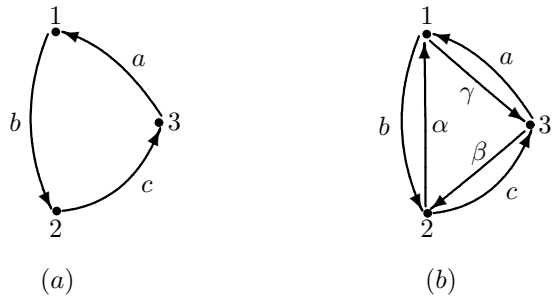

\vspace*{0mm}
\begin{center}
\setlength{\unitlength}{0.85mm}
\hspace{0mm}
\input cheb_5.lp
\caption{Weighted digraphs on three vertices: (a) a weighted cycle $C_3$; (b) a complete digraph $K_3$.\label{f_3Wei}}
\vspace{-1.0em}
\end{center}
\end{figure}

First, consider an arbitrary weighted directed cycle $C_3$ on three vertices, see Fig.~\ref{f_3Wei}(a).
Recall that by Theorem~\ref{t_G1G2} the \emph{unweighted\/} directed cycle is essentially cyclic.

\begin{proposition}
\label{t_3cycle}
The weighted $3$-cycle $C_3$ 
is essentially cyclic if and only if the square roots $\sqrt{a},$ $\sqrt{b},$ and $\sqrt{c}$
of its arc weights satisfy the strict triangle inequality$,$ namely$:$ 
\[
\sqrt{a}<\sqrt{b}+\sqrt{c},\quad \sqrt{b}<\sqrt{a}+\sqrt{c},\,\;\;\text{and}\/\;\;\;\sqrt{c}<\sqrt{a}+\sqrt{b}.
\]
\end{proposition}

Proposition~\ref{t_3cycle} follows from Theorem~\ref{t_3di} below. It
can be interpreted as follows: in an essentially cyclic digraph $C_3$, the weight of any arc
is not large enough to ``overpower'' the remaining part of the cycle. Or, in more precise terms, $C_3$ is essentially
cyclic whenever the square roots of its arc weights are the lengths of the sides of a non-degenerate triangle.

Now consider the complete weighted digraph (without loops) on three vertices, $K_3$ (Fig.~\ref{f_3Wei}(b)).
\bigskip

The Laplacian characteristic equation of this digraph,
\begin{equation*}
    \det
    \begin{pmatrix}
\la-b-\g &        b &       \g\\
      \a & \la-c-\a &        c\\
       a &       \b & \la-a-\b\\
    \end{pmatrix}
=0,
\end{equation*}
reduces to
\begin{equation*}
    \la(\la^2-(a+b+c+\a+\b+\g)\la+(ab+bc+ca+\a\b+\b\g+\g\a+a\a+b\b+c\g))=0.
\end{equation*}
The digraph is essentially cyclic if and only if $D<0,$ where
\begin{equation*}
   D=(a+b+c+\a+\b+\g)^2-4(ab+bc+ca+\a\b+\b\g+\g\a+a\a+b\b+c\g).
\end{equation*}

Equivalently,
\[
D=(a-\a)^2-2(a-\a)(b+c-\b-\g)+(b-c-\b+\g)^2.
\]
This quadratic trinomial in $a-\a$ is negative iff its roots are real and $a-\a$ lies strictly between them.
The roots
\begin{eqnarray}
\nonumber
(a-\a)_{1,2} &=& (b+c-\b-\g)\pm\sqrt{(b+c-\b-\g)^2-(b-c-\b+\g)^2}\\
             &=& (b-\b)+(c-\g)\pm2\sqrt{(b-\b)(c-\g)}
\label{roots3}
\end{eqnarray}
are real and unequal iff
\beq
\label{posir}
(b-\b)(c-\g)>0.
\eeq

Assuming that \eqref{posir} is satisfied, first consider the case of $b>\b$ and $c>\g$. In this case, \eqref{roots3}
reduces to the entirely real expression
\[
(a-\a)_{1,2}=(\sqrt{b-\b}\pm\sqrt{c-\g})^2.
\]
Then the inequality $D<0$, i.\,e. the essential cyclicity of $K_3$, amounts to
\[
\left|\sqrt{b-\b}-\sqrt{c-\g}\right|<\sqrt{a-\a}<\sqrt{b-\b}+\sqrt{c-\g},
\]
which is the strict triangle inequality for $\sqrt{a-\a}$, $\sqrt{b-\b},$ and $\sqrt{c-\g}$.

The case of $b<\b$ and $c<\g$ is considered similarly; as a result we obtain the following theorem.

\begin{theorem}
\label{t_3di}
Let the matrix of arc weights of a weighted digraph\/ $\G$ be
$\,W=
\left(
\begin{array}{rrr}
     0&     b&     \g\\
    \a&     0&      c\\
     a&    \b&      0\\
\end{array}
\right)
$.
Then $\G$ is essentially cyclic if and only if either\\
{\rm (i)\phantom{i}}
       $\sqrt{a-\alpha},$ $\sqrt{b-\beta},$ and $\sqrt{c-\gamma}$ are real and satisfy the strict triangle inequality or\\
{\rm (ii)}
       $\sqrt{\alpha-a},$ $\sqrt{\beta-b},$ and $\sqrt{\gamma-c}$ are real and satisfy the strict triangle inequality.
\end{theorem}

This criterion corresponds to a certain intuitive sense of essential cyclicity.

According to Theorem~\ref{t_3di},
for a complete weighted digraph on three vertices, $K_3$, to be essentially cyclic, it is necessary that it obeys\x 
the strict triangle inequality applied to certain values rather attached to the vertices than to the pairs of them.
Indeed, such a value is the square root of the weight difference for the two arcs converging to \emph{the same\/} vertex.

The problem of characterizing the essentially cyclic weighted digraphs becomes much more difficult with
the increase of the number of vertices. However, some of the corresponding conditions involve the triangle inequality for
the roots of the arc weights as well.

\begin{figure}[htb]
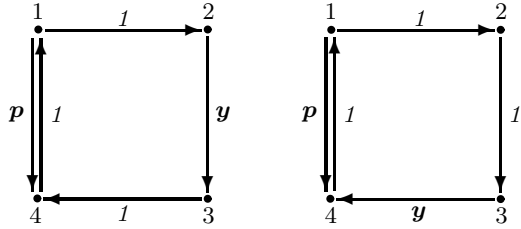

\vspace*{0mm}
\begin{center}
\setlength{\unitlength}{0.75mm}
\hspace{0mm}
\input cheb_6a.lp
\vspace{-1em}
\caption{Two cospectral weighted digraphs with ring structure on four vertices.\label{f_4Wei}}
\vspace{-1.0em}
\end{center}
\end{figure}

Consider the first weighted digraph in Fig.~\ref{f_4Wei}.
Note that the corresponding unweighted digraph is essentially cyclic by Theorem~\ref{t_gener}.
The Laplacian characteristic equation for the weighted digraph,
\begin{equation*}
    \det
    \begin{pmatrix}
\la-p-1 &    1 &    0 &    p\\
      0 &\la-y &    y &    0\\
      0 &    0 &\la-1 &    1\\
      1 &    0 &    0 &\la-1\\
    \end{pmatrix}
=0,
\end{equation*}
as well as that for the second weighted digraph shown in Fig.~\ref{f_4Wei}, reduces to the form
\begin{equation*}
\la(\la^3-(y+q)\la^2+(qy+q)\la-(qy+1))=0, \;\;\mbox{where}\;\; q=p+3.
\end{equation*}

The digraph is essentially cyclic whenever $D<0,$ where
\beq
\label{e_det3}
D=\frac{4(b^2-3c)^3-(2b^3-9bc+27d)}{27},
\eeq
\begin{equation*}
b=-(y+q),\quad c=qy+q,\;\;\;\text{and}\;\;\;d=-(qy+1).
\end{equation*}

Substituting these expressions for $b$, $c$, and $d$ in \eqref{e_det3} we obtain
\beq
\label{D4}
D=q(q-4)y^4-(2q^3-8q^2+4)y^3+q(q^3-2q^2-8q+6)y^2-2q(q+2)(q-3)^2y+(q+1)(q-3)^3.
\eeq

To solve the inequality $D<0$ w.r.t. $y$, one can find the real roots of the polynomial~\eqref{D4}
treating $q$ as a parameter. However, the expressions of these roots as functions in $q$ are rather cumbersome,
as well as the inverse representations of $q$ via~$y$.

For example, when $p=3$, the roots of the equation $D=0$ are:
\[
y_{1,2}=\left(37-Q\pm\sqrt{290-36z-504/z+3454/Q}\right)\!/12,
\]
where
\[
Q=\sqrt{36z+145+504/z},\quad z=0.5\sqrt[3]{671+65\sqrt{65}}.
\]
Therefore it does not seem to be easy to formulate a simple criterion (such as Theorem~\ref{t_3di}) of essential cyclicity
for the digraphs of this kind.
Solving the problem numerically, we obtain that for
$p=3$, the digraph is essentially cyclic, i.\,e.\ $D<0$, when $0.266<y<2.441$ (approximately). So the essential cyclicity can be
suppressed by either increase or decrease of the arc weight~$y$.

\begin{figure}[htb]
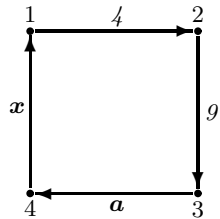

\vspace*{0mm}
\begin{center}
\setlength{\unitlength}{0.72mm}
\noindent
\hspace{-2em}
\input cheb_7a.lp
\vspace{-1em}
\caption{A cycle on four vertices with two variable arc weights, $a$ and $x$.\label{f_4Wei1}}
\vspace{-1.0em}
\end{center}
\end{figure}

Finally, consider a weighted cycle on four vertices with two variable weights (Fig.~\ref{f_4Wei1}).

The corresponding Laplacian characteristic polynomial is
\[
f(\la)=\la(\la^3-(13+x+a)\la^2+(36+13x+13a+ax)\la-36x-36a-13ax).
\]

The boundary of the domain on the $(a,x)$ plane corresponding to the essentially cyclic digraphs 
is specified by the equation
\begin{eqnarray*}
&-&a^2x^2(x-a)^2+26(x+a)\bigl(ax(x-a)^2+25(x^2+a^2)+58ax+900\bigr)\\
&+&870a^2x^2-241(x^2+a^2)(2ax+25)-25(x^4+a^4)-3934ax-32400=0.
\end{eqnarray*}
\begin{figure}[htb]
\hspace{7em}\includegraphics[width=4.168in,height=4.2in]{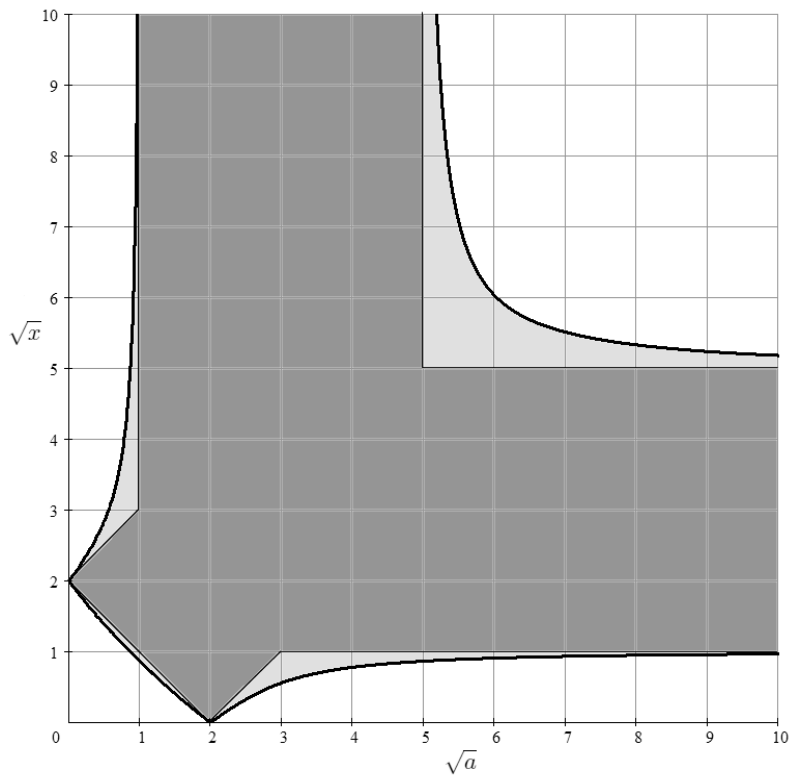} \vspace{-1.5em}
\caption{For the weighted digraph shown in Fig.~\ref{f_4Wei1}, the domain where the square roots of the three smaller arc weights satisfy the triangle inequality is filled in dark grey; the domain corresponding to the essentially cyclic digraphs is the union of the above ``triangle inequality domain'' and the four fringing domains filled in light grey.} \label{f_1ort}
\end{figure}
The polynomial on the left-hand side is not the product of polynomials with rational coefficients and smaller degrees.
The solutions $(a,x)$ of the above equation in the non-negative quadrant are plotted in Fig.~\ref{f_1ort} with the ``rooted'' scales $\sqrt{a}$ and $\sqrt{x}$.

The domain corresponding to the essentially cyclic digraphs is filled.
The subdomain filled in dark grey is the locus of points $(\sqrt{a},\sqrt{x})$ for which the square roots of
the three smallest arc weights satisfy the triangle inequality.
The locus of points $(\sqrt{a},\sqrt{x})$ that correspond to the essentially cyclic digraphs is wider: it also contains the four subdomains filled in light grey. 

Thus, for this cyclic digraph, the triangle inequality for the square roots of the three smaller arc weights is a
sufficient, but not necessary condition of essential cyclicity.
In other words, the required criterion of essential cyclicity is a kind of {\em relaxed\/} triangle inequality.
This relaxed inequality turns into the triangle inequality as the fourth (largest) arc weight tends to infinity
or the smallest weight tends to zero; the relaxation is maximal when the largest weight becomes equal to the second
largest weight. It can be conjectured 
that the triangle inequality for the square roots of the three smallest arc weighs is a sufficient condition of essential cyclicity for the whole
class of weighted 4-cycles.

As one can see, even for the weighted digraphs on four vertices, the problem of characterizing essential cyclicity in terms
of graph topology is non-trivial.

\section{Concluding remarks}
\label{s_concl}

We conclude with two side 
remarks.

1. According to the matrix tree theorem (see, e.g., Theorem~16.9 in \cite{Harary}), for every digraph $\G$,
the cofactor of each entry in the $i$th row of the Laplacian matrix is equal to the number of spanning
converging trees (also called in-arborescences) rooted at vertex~$i$. The total number of in-arborescences is
equal to the sum of the cofactors in any column of~$L$.

Thus, the matrix tree theorem provides a general approach to computing the number of in-arborescences in a
digraph and the number of spanning trees in a graph. For certain classes of graphs, this approach leads to
explicit formulas which can be obtained using the Chebyshev polynomials of the second kind. In
\cite{BoeschProdinger}, this method was applied to the wheels, fans, M\"obius ladders, etc. In
\cite{ZhangYongGolin}, the Chebyshev polynomials were used for finding the number of spanning trees for
certain classes of graphs including circulant graphs with fixed and non-fixed jumps. The results of the
present paper can be used to obtain representations for the number of converging trees in the digraphs with
ring structure.

Suppose that $t_i^n$ is the number of spanning converging trees in the digraph whose Laplacian matrix is
$L''_n$ (Eq.~\eqref{L''}). Then summing up the cofactors of the last column of $L''_n$ and having in mind that
(i)~$\det M_n=(-1)^nZ_n(0)=(-1)^nP_{2n}(0)=1$,\, 
(ii)~$P_n(x)=\det(xI-C_n)$, where
\[
C_n=\begin{pmatrix}
 0    & 1    &      &      &      &\0L\\
 1    & 0    & 1    &      &      &   \\
      &\ddots&\ddots&\ddots&      &   \\
      &      &\ddots&\ddots&\ddots&   \\
      &      &      & 1    & 0    & 1 \\
\0L   &      &      &      & 1    & 0\\
    \end{pmatrix}
\]
(see, e.\,g., \cite[Theorem~1.11]{Paszkowski}), and\,
(iii)~$P_k(2)=k+1$ (which equals the limit from the left at $x=2$ of the expression~\eqref{080207eq1}) we obtain the
following representation for $t_i^n$:
\beq
\label{250509eq1}
t_i^n
=\suml_{k=0}^{i-1}P_k(2)+\suml_{k=0}^{n-i-1}P_k(2)
={1\over 2}\bigl(i^2+n+(n-i)^2\bigr).
\eeq

Since $t_i^n$ is the product of the eigenvalues of $L''_n$, except for a zero eigenvalue,
for the non-trivial digraphs with ring structure that have completely real spectra,
Eq.~\eqref{250509eq1} and Theorem~\ref{280207th2} provide the following expressions for $t_i^n$ and simultaneously
trigonometric identities:
\begin{eqnarray*}
&&t_{n/2}^n
=\left(
\prod_{k=1}^{  n/2-1}\!\!2\cos{\pi k\over n  }\,\,
\prod_{k=1}^{  n/2  }    2\cos{\pi k\over n+2}
\right)^{\!2}
={n(n+2) \over 4},\;\text{if}\;n\;\text{is even};\\
&&t_{(n-1)/2}^n=t_{{n+1}/2}^n
=\left(
\prod_{k=1}^{(n-1)/2}\!\!2\cos{\pi k\over n+1}
\right)^{\!4}
={(n+1)^2\over 4},\;\text{if}\;n\;\text{is odd}.
\end{eqnarray*}

2. Let $L_n^c$ be the Laplacian matrix of the undirected cycle on $n$ vertices. By suitable indexing of the vertices,
$L_n^c$ can be presented as a matrix different from $M_n$ \eqref{280609eq3} in the $(1,1)$ entry, which in
$L_n^c$ is~$1$.

Expanding $\det(\la I-L_n^c)$ along the first row and using \eqref{010207eq1}, Lemma~\ref{010207prop1},
and \eqref{e_P2m1}, for all $\la\in\,]0,\,4]$ one has:
\begin{eqnarray*}
\nonumber
\Delta_{L^c_n}(\la)
&=&(\la-1)Z_{n-1}(\la) Z_{n-2}(\la)
 =        Z_n    (\la)+Z_{n-1}(\la)
 =P_{2n}(\sqrt{\la})+P_{2(n-1)}(\sqrt{\la})\\
&=&      \sqrt{\la}\,P_{2 n-1} (\sqrt{\la})
 =\la\prod_{k=1}^{ n-1}\Big(\la-4\cos^2{\pi k\over 2n}\Big)
 =   \prod_{k=1}^{ n  }\Big(\la-4\cos^2{\pi k\over 2n}\Big).
\end{eqnarray*}

Thus, the roots of $\Delta_{L^c_n}$ are $(4\cos^2{\pi k\over 2n},\; k=\1n)$. Obviously,
$(4\sin^2{\pi k\over 2n},\;k=0\cdc n-1)$ is a different representation of the same spectrum.
The latter representation was taken as ``ready-made'' and then proved in \cite{Fiedler72} and~\cite{AndersonMorley}.
The above reduction to Chebyshev polynomials provides a derivation of this result.
Another derivation can be obtained using item~1 of Theorem~\ref{t_G1G2} and the representation
$C=\vec{C}\cup\vec{C}^{n-1}$, where $\vec{C}$ is the directed cycle and $\vec{C}^{n-1}$ its $(n-1)$st power. On connections of the
adjacency characteristic polynomials with Chebyshev polynomials, see \cite[\S\,{2.6}]{CvetkovicDoob80}.

\section*{Acknowledgements}

This work was partially supported by RFBR Grant 09-07-00371 and the RAS Program ``Development of Network and Logical
Control in Conflict and Cooperative Environments.''


\end{document}